\documentclass{amsart}

\def\Q{{\mathbf Q}}
\def\Z{{\mathbf Z}}
\def\C{{\mathbf C}}

\def\r{{\mathfrak G}}
\def\A{{\mathbf A}}
\def\Gal{\mathrm{Gal}}
\def\End{\mathrm{End}}
\def\Aut{\mathrm{Aut}}
\def\Hom{\mathrm{Hom}}
\def\Lie{\mathrm{Lie}}

\def\ppr{\mathrm{prime}}

\def\fchar{\mathrm{char}}
\def\GL{\mathrm{GL}}
\def\Im{\mathrm{Im}}
\def\dim{\mathrm{dim}}

\def\O{{\mathcal O}}
\def\X{{A}}
\def\Y{{B}}
\def\CC{{U}}
\def\k{k}
\def\kk{K}

\def\P{{\mathcal P}}

\def\bigtilde{\widetilde}
\def\bmu{\hbox{\rm\rlap {$\mu$}\kern-.07em {$\mu$}}}

\def\downdiag#1#2#3#4#5#6{
   \begin{center}
   \begin{picture}(85,50)
   \put(-3,40){\makebox(0,0){$#1$}}
   \put(80,40){\makebox(0,0){$#2$}}
   \put(80,5){\makebox(0,0){$#3$}}
   \put(18,40){\vector(1,0){38}}
   \put(80,30){\vector(0,-1){15}}
   \put(18,33){\vector(2,-1){41}}
   \put(37,44){\makebox(0,0)[b]{\scriptsize $#4$}}
   \put(36,18){\makebox(0,0)[r]{\scriptsize $#5$}}
   \put(83,24){\makebox(0,0)[l]{\scriptsize $#6$}}
   \end{picture}
   \end{center}}

\newtheorem{thm}{Theorem}[section]
\newtheorem{lem}[thm]{Lemma}
\newtheorem{cor}[thm]{Corollary}
\newtheorem{prop}[thm]{Proposition}
\newtheorem{claim}{Lemma}[subsection]
\theoremstyle{definition}
\newtheorem{defn}[thm]{Definition}

\newtheorem{rem}[thm]{Remark}
\newtheorem{remss}[claim]{Remark}

\title[Connectedness extensions for abelian varieties]
{Connectedness extensions for abelian varieties}

\author[A.\ Silverberg]{A.\ Silverberg*}
\thanks{* Partially supported by the National
Science Foundation.}
\address{Department of Mathematics, Ohio State University,
231 W.\ 18 Avenue,
Columbus, Ohio 43210--1174, USA}
\email{silver\char`\@math.ohio-state.edu}

\author[Yu. G. Zarhin]{Yu. G. Zarhin**}
\thanks{** Partially supported by the Ohio State
University Mathematics Research Institute and the Deutsche
Forschungsgemeinschaft.}
\address{Department of Mathematics, Pennsylvania State University,
University Park, PA 16802, USA,
\newline
\indent Institute for Mathematical Problems in Biology,
Russian Academy of Sciences, Push\-chino, Moscow Region, 142292, Russia}
\email{zarhin\char`\@math.psu.edu}

\begin{document}

\maketitle

\section{Introduction}

Suppose $\X$ is an abelian variety defined over a field $F$,
$\ell$ is a  prime number, and $\ell \neq \fchar(F)$.
Let $F^s$ denote a separable closure of $F$, let
$T_\ell(\X) = {\displaystyle \lim_\leftarrow \X_{\ell^r}}$
(the Tate module), let
$V_\ell(\X) = T_\ell(\X) \otimes_{\Z_\ell}\Q_\ell$, and let
$\rho_{\X,\ell}$
denote the $\ell$-adic representation
$$\rho_{\X,\ell} : \Gal(F^s/F) \to \Aut(T_\ell(\X)) \subseteq
\Aut(V_\ell(\X)).$$
If $L$ is an extension of $F$ in $F^s$, let
$G_{L,\X}$ denote the image of $\Gal(F^s/L)$ under $\rho_{\X,\ell}$.
Let $\r_\ell(F,\X)$
denote the algebraic envelope of the image of $\rho_{\X,\ell}$,
i.e., the Zariski closure of $G_{F,\X}$ in
$\Aut(V_\ell(\X)) \cong \GL_{2d}(\Q_\ell)$,
where $d = \dim(\X)$.

Let $F_{\Phi,\ell}(\X)$ be the smallest extension $F'$ of $F$ such that
$\r_\ell(F',\X)$ is connected. We call this extension the
$\ell$-connectedness extension, or connectedness extension.

The algebraic group $\r_\ell(F,\X)$
and the field $F_{\Phi,\ell}(\X)$ were
introduced by Serre (\cite{serre},
\cite{resume}, \cite{serremotives}), who
proved that if $F$ is a global field or a finitely
generated extension of $\Q$,
then $F_{\Phi,\ell}(\X)$ is independent of $\ell$
(see also \cite{LPinv}, \cite{LPmathann}, \cite{LP}). In such cases,
we will denote the field $F_{\Phi,\ell}(\X)$ by $F_\Phi(\X)$. For
every integer $n \ge 3$ we have
$$F_{\Phi}(\X) \subseteq F(\X_n)$$
(see \cite{Borovoi}, \cite{Borovoisb}, Proposition 3.6 of
\cite{Chi}, and \cite{conn}). Larsen and Pink \cite{LP}
recently proved that for every integer $n \ge 3$,
$$F_\Phi(\X) = \bigcap_{\ppr~~p \geq n} F(\X_p).$$
In \cite{conn} we found conditions for the connectedness
of $\r_\ell(F,\X)$, while in \cite{connexts} we used connectedness
extensions and Serre's $\ell$-in\-de\-pen\-dence results
to obtain $\ell$-in\-de\-pen\-dence results for the intersection of
$\r_\ell(F,\X)(\Q_\ell)$ with the torsion subgroup of the center of
$\End(\X) \otimes \Q$.

Let $F(\End(\X))$
denote the smallest extension of $F$ over which all the endomorphisms
of $\X$ are defined. Then
(see Proposition 2.10 of \cite{conn}),
$$F(\End(\X)) \subseteq F_{\Phi,\ell}(\X).$$
Therefore, $\r_\ell(F,\X)$ fails to be connected when the ground field is
not a field of definition for the endomorphisms of $\X$.
For example, if $F$ is a subfield of $\C$, and
$\X$ is an elliptic curve over $F$ with complex multiplication
by an imaginary quadratic field $K$ which is not contained in $F$, then
$F \neq KF = F(\End(\X)) \subseteq F_{\Phi,\ell}(\X)$.
More generally, if $\X$ is an abelian variety
of CM-type, and ${\tilde K}$ is the
reflex CM-field, then $F(\End(\X)) \supseteq {\tilde K}$; if
${\tilde K}$ is not contained in $F$ then
$F \neq F(\End(\X)) \subseteq F_{\Phi,\ell}(\X)$.
It is therefore natural to enlarge
the ground field $F$ so that it is a field of definition for the
endomorphisms of $\X$.

By enlarging the ground field, we may assume that
$F = F(\End(\X)) = F_{\Phi,\ell}(\X)$. We then consider the $F$-forms
$\Y$ of $\X$ such that $F = F(\End(\Y))$. For such $\Y$, we describe
the connectedness extensions
$F_{\Phi,\ell}(\Y)/F$ (see \S\ref{exclasses}, especially
Theorem \ref{discond} and Corollary \ref{discondcor}).
Properties of Mumford-Tate groups given in \S\ref{notation} allow
us to obtain explicit information about the connectedness extensions
$F_{\Phi,\ell}(\Y)/F$ under additional conditions (see
Theorems \ref{twistX} and \ref{twistE}).
Our conditions in Theorems \ref{twistX} and \ref{twistE}
are based on Weil's philosophy in \cite{WeilHodge} whereby
exceptional Hodge classes arise from certain abelian varieties
that have a CM-field embedded in their endomorphism algebras.
In \S \ref{examples} we use the results of \S\ref{exclasses}
to explicitly compute non-trivial connectedness extensions in
special cases.

\noindent {\bf Acknowledgments:}
The authors would like to thank the Mathematische Institut der
Universit\"at Erlangen-N\"urnberg for its hospitality.

\section{Definitions, notation, and lemmas}
\label{notation}

Let $\Z$, $\Q$, and $\C$ denote respectively the integers, rational
numbers, and complex numbers. If $r$ is an integer, then
$\Q(r)$ denotes the rational Hodge structure of weight $-2r$ on $\Q$ (see
\S 1 of \cite{sln900}). If $a$ and $b$ are integers, let $(a,b)$
denote the greatest common divisor of $a$ and $b$.
If $F$ is a field, let $F^s$ denote a separable closure of $F$ and let
${\bar F}$ denote an algebraic closure of $F$. If $\X$ is an abelian variety
over a field $F$, write $\End_F(\X)$
for the set of endomorphisms of $\X$ which are defined over $F$, let
$\End(\X) = \End_{F^s}(\X)$, and let $\End^0(\X) = \End(\X) \otimes_\Z \Q$.
Let $Z_\X$ denote the center of $\End(\X)$.
If $G$ is an algebraic
group, let $G^0$ denote the identity connected component.

\begin{lem}[Lemma 2.7 of \cite{conn}]
\label{conncomp}
If $\X$ is an abelian variety over a field $F$, $L$ is a finite
extension of $F$ in $F^s$, and $\ell$ is a prime number, then
$$\r_\ell(L,\X) \subseteq \r_\ell(F,\X)
\text{ and } \r_\ell(L,\X)^0 = \r_\ell(F,\X)^0.$$
In particular, if $\r_\ell(F,\X)$
is connected, then $\r_\ell(F,\X) = \r_\ell(L,\X)$.
\end{lem}

\begin{lem}
\label{conncomplem}
Suppose $\X$ and $\Y$ are abelian varieties over a field $F$, $L$ is a finite
extension of $F$ in $F^s$, $\ell$ is a prime number,
$\ell \neq \fchar(F)$, $\r_\ell(F,\X)$ is
connected, and $\X$ and $\Y$ are isomorphic over $L$. Then:
\begin{enumerate}
\item[{(i)}] $\r_\ell(F,\Y)^0 = \r_\ell(F,\X)$, and
\item[{(ii)}]$\r_\ell(L,\Y)$ is connected, i.e.,
$F_{\Phi,\ell}(\Y) \subseteq L$.
\end{enumerate}
\end{lem}

\begin{proof}
Since $\X$ and $\Y$ are isomorphic over $L$, and $\r_\ell(F,\X)$ is
connected, we have
$$\r_\ell(L,\Y) = \r_\ell(L,\X) = \r_\ell(F,\X) = \r_\ell(F,\X)^0$$
$$= \r_\ell(L,\X)^0 = \r_\ell(L,\Y)^0 = \r_\ell(F,\Y)^0,$$
using Lemma \ref{conncomp}. The result follows.
\end{proof}

\begin{prop}
\label{endaut}
Suppose $\X$ and $\Y$ are abelian varieties over a field $F$, $L$ is a field
extension of $F$ in $F^s$, and $f: \X \to \Y$ is an isomorphism defined over
$L$.  Suppose that for every $\sigma \in \Gal(F^s/F)$, the element
$f^{-1}\sigma(f)$  of $\Aut(\X)$ commutes with
every element of $\End_L(\X)$. Then $\End_F(\X) \cong \End_F(\Y)$.
\end{prop}

\begin{proof}
Define an isomorphism $\varphi : \End_L(\X) \to \End_L(\Y)$ by
$\varphi(\beta) = f\beta f^{-1}$. For every $\beta \in \End_L(\X)$ and $\sigma
\in \Gal(F^s/F)$, we have
$f^{-1}\sigma (f)\beta = \beta f^{-1}\sigma(f)$.
Therefore, $\sigma(f\sigma^{-1}(\beta) f^{-1}) = f\beta f^{-1}$.
Thus, $\beta \in \End_F(\X)$ if and only if $f\beta f^{-1} \in \End_F(\Y)$.
In other words, the restriction of $\varphi$ to $\End_F(\X)$ induces an
isomorphism onto $\End_F(\Y)$.
\end{proof}

As a corollary we have the following result. See also Lemma 5.1 of
\cite{szisog}.

\begin{cor}
\label{cocy}
Suppose $\X$ is an abelian variety over a field $F$. If an element
of $H^1(\Gal(F^s/F),\Aut(\X))$ is represented by a cocycle $c$
with values
in the center of $\End^0(\X)$, and $\Y$ is the twist of $\X$ by $c$,
then $\End_F(\X) \cong \End_F(\Y)$.
\end{cor}

\begin{proof}
The cocycle $c$ defines an isomorphism $f : \X \to \Y$ such that for every
$\sigma \in \Gal(F^s/F)$, $f^{-1}\sigma(f) = c(\sigma)$. We apply
Proposition \ref{endaut}.
\end{proof}

\begin{lem}
\label{cocychar}
Suppose $\X$ is an abelian variety over a field $F$,
$c$ is a cocycle on $\Gal(F^s/F)$ with values in $\Aut(\X)$,
$\Y$ is the twist of $\X$ by $c$,
and $F = F(\End(\X)) = F(\End(\Y))$.
Then $c$ is a character with values in $Z_\X^\times$,
where $Z_\X$ denotes the center of $\End(\X)$.
\end{lem}

\begin{proof}
Since $\Gal(F^s/F)$ acts trivially on $\End(\X)$,
the cocycle $c$ is a homomorphism. Let $f : \X \to \Y$ be the
isomorphism induced by $c$. Then $c(\sigma) = f^{-1}\sigma(f)$
for every $\sigma \in \Gal(F^s/F)$. Since
$F = F(\End(\X)) = F(\End(\Y))$, it easily follows that
$c(\sigma) \in Z_\X$ and $c(\sigma)^{-1} \in Z_\X$.
\end{proof}

\begin{rem}
\label{twistyrem}
If an abelian variety $\Y$ over $F$ is the twist of an abelian variety
$\X$ by $c \in H^1(\Gal(F^s/F),\Aut(\X))$ then one may easily check that
the Galois module $\Y(F^s)$ is the twist by $c$ of the Galois module
$\X(F^s)$, and therefore the Galois module $V_\ell(\Y)$ is the twist by
$c$ of the Galois module $V_\ell(\X)$.
\end{rem}

We define the Mumford-Tate group of a complex abelian variety $\X$
(see \S2 of \cite{Ribet} or \S6 of \cite{Izv}). If $\X$ is a complex abelian
variety, let $V = H_1(\X(\C),\Q)$ and consider the Hodge decomposition
$V \otimes \C = H_1(\X(\C),\C) = H^{-1,0} \oplus H^{0,-1}$.
Define a homomorphism $\mu : {\mathbf G}_m \to GL(V)$ as follows. For
$z \in \C$, let $\mu(z)$ be the automorphism of $V \otimes \C$ which is
multiplication by $z$ on $H^{-1,0}$ and is the identity on $H^{0,-1}$.

\begin{defn}
The {\em Mumford-Tate group} $MT_\X$ of $\X$ is the smallest
algebraic subgroup of $GL(V)$, defined over $\Q$, which after extension of
scalars to $\C$ contains the image of $\mu$.
\end{defn}

It follows from the definition that $MT_\X$ is connected.

Define a homomorphism $\varphi : {\mathbf G}_m \times {\mathbf G}_m \to GL(V)$
as  follows. For $z, w \in \C$, let $\varphi(z,w)$ be the automorphism of
$V \otimes \C$ which is multiplication by $z$ on $H^{-1,0}$ and is
multiplication by $w$ on $H^{0,-1}$. Then $MT_\X$ can also be defined as
the smallest algebraic subgroup of $GL(V)$, defined over $\Q$, which after
extension of scalars to $\C$ contains the image of $\varphi$. The equivalence
of the definitions follows easily from the fact that $H^{-1,0}$ is the complex
conjugate of $H^{0,-1}$. (See \S 3 of \cite{serrereps}, where $MT_\X$ is called
the Hodge group. See also \S 6 of \cite{Izv}.)

If $\X$ is an abelian variety over a subfield $F$ of $\C$, we fix an
embedding of ${\bar F}$ in $\C$. This gives an identification of
$V_\ell(\X)$ with $H_1(\X,\Q)\otimes\Q_\ell$, and allows us to view
$MT_\X \times \Q_\ell$ as a linear $\Q_\ell$-algebraic subgroup of
$GL(V_\ell(\X))$. Let
$$MT_{\X,\ell} = MT_\X \times_\Q \Q_\ell.$$
Then $MT_{\X}(\Q_\ell) = MT_{\X,\ell}(\Q_\ell)$.
The Mumford-Tate conjecture for abelian varieties (see \cite{serrereps})
may be reformulated as the equality of $\Q_\ell$-algebraic groups,
$\r_\ell(F,\X)^0 = MT_{\X,\ell}$.

\begin{thm}[Piatetski-Shapiro \cite{ps}, Deligne \cite{sln900}, Borovoi
\cite{Borovoisb}]
\label{psdcor}
If $\X$ is an abelian variety over a finitely generated extension $F$ of
$\Q$, then $\r_\ell(F,\X)^0 \subseteq MT_{\X,\ell}$.
\end{thm}

In \S\ref{exclasses} it will be helpful to use a slightly
different version of the Mumford-Tate group, as defined by Deligne
(see p.~43 and pp.~62--63 of \cite{sln900}). We will denote this group
${\bigtilde{MT}}_\X$. (See also pp.~466--467 of \cite{Milne} for a
comparison between $MT_\X$ and ${\bigtilde{MT}}_\X$.)
Letting $V^*$ be the dual of $V$, then
$T = V^{\otimes p} \otimes (V^*)^{\otimes q} \otimes \Q(r)$ has a Hodge
structure of weight $q - p - 2r$. If $\nu \in {\mathbf G}_m$, let $\nu$ act on
$\Q(1)$ as $\nu^{-1}$, and we obtain a canonical action of
$GL(V) \times {\mathbf G}_m$ on $T$. (Note that $V^* \cong V \otimes \Q(1)$,
since $V$ is a polarized Hodge structure of weight $-1$.)

\begin{defn}
The group ${\bigtilde{MT}}_\X$ is the subgroup of
$GL(V) \times {\mathbf G}_m$ consisting of the elements which fix all
rational tensors of bidegree $(0,0)$ belonging to any $T$.
\end{defn}

\begin{lem}[Proposition 3.4 of \cite{sln900}]
\label{murem}
The algebraic group ${\bigtilde{MT}}_\X$ is the smallest algebraic subgroup of
$GL(V) \times {\mathbf G}_m$ defined over $\Q$ which, after extension of
scalars to $\C$, contains the image of
$(\mu,{\mathrm{id}}) : {\mathbf G}_m \to GL(V) \times {\mathbf G}_m$.
\end{lem}

If $F$ is a field and $\ell$ is a prime
number different from $\fchar(F)$, let
$$\chi_\ell : \Gal(F^s/F) \to \Z_\ell^\times \subset \Q_\ell^\times$$
denote the cyclotomic character. If $r$ is an integer, then the
$\Gal(F^s/F)$-module $\Q_\ell(r)$ is the $\Q_\ell$-vector space $\Q_\ell$
with Galois action defined by the character $\chi_\ell^r$. We have
$\Q_\ell(r) = \Q(r) \otimes_\Q \Q_\ell$ (see \S 1 of \cite{sln900}).
Suppose $\X$ is an abelian variety over $F$. Let $V_\ell = V_\ell(\X)$ and
let $V_\ell^*$ be the dual of $V_\ell$.
If $\nu \in {\mathbf G}_m$, let $\nu$ act on
$\Q_\ell(1)$ as $\nu^{-1}$. We obtain a canonical action of
$GL(V_\ell) \times {\mathbf G}_m$ on
$V_\ell^{\otimes p} \otimes (V_\ell^*)^{\otimes q} \otimes \Q_\ell(r)$.
Define
$${\tilde \rho}_{\X,\ell} : \Gal(F^s/F) \to \Aut(V_\ell) \times \Q_\ell^\times
= \Aut(V_\ell) \times {\mathbf G}_m(\Q_\ell)$$
by ${\tilde \rho}_{\X,\ell}(\sigma) =
(\rho_{\X,\ell}(\sigma),\chi_\ell^{-1}(\sigma))$.

\begin{defn}
Let ${\tilde \r}_\ell(F,\X)$ denote the smallest
$\Q_\ell$-algebraic subgroup of
$$GL(V_\ell) \times {\mathbf G}_m$$ whose group of $\Q_\ell$-points contains
the image of ${\tilde \rho}_{\X,\ell}$.
\end{defn}

If $\X$ is a complex abelian variety, then a polarization on $\X$ (i.e.,
the imaginary part of a Riemann form) produces an element $E$ of
$\Hom(\wedge^2V,\Q(1))$ which is a rational tensor of bidegree $(0,0)$.
If $\X$ is an abelian variety over an arbitrary field $F$, then a
polarization on $\X$ defined over $F$ defines a $\Gal(F^s/F)$-invariant
element $E_\ell$ of $\Hom(\wedge^2V_\ell,\Q_\ell(1))$
(since the Weil pairing is $\Gal(F^s/F)$-equivariant).
If ${\bar F}$ is a
subfield of $\C$, and we fix a polarization on $\X$ defined over $F$,
then the line generated by $E_\ell$ in $\Hom(\wedge^2V_\ell,\Q_\ell(1))$
is the extension of scalars to $\Q_\ell$ of the line generated by $E$
in $\Hom(\wedge^2V,\Q(1))$.
(See p.~237 of \cite{MumfordAV}, especially the last sentence.)

The following result implies that the projection map
$GL(V) \times {\mathbf G}_m \to GL(V)$ induces an isomorphism from
${\bigtilde{MT}}_\X$ onto $MT_\X$. Since we were not able to find a proof
in the literature, we have included one for the benefit of the
reader.

\begin{prop}
\label{MTtilde}
If $\X$ is a complex abelian variety, then there exists a (unique)
character $\gamma : MT_\X \to {\mathbf G}_m$ such that
${\bigtilde{MT}}_\X$ is the graph of $\gamma$.
\end{prop}

\begin{proof}
Let $p_1$ and $p_2$ denote the
projection maps from $GL(V) \times {\mathbf G}_m$ onto $GL(V)$ and ${\mathbf
G}_m$,
respectively. By Lemma \ref{murem}, $MT_\X$ is the image of
${\bigtilde{MT}}_\X$
under $p_1$.
Fix a polarization on $\X$. The polarization generates a line $D$ in the
$\Q$-vector space $\Hom(\wedge^2V,\Q(1))$, on which ${\bigtilde{MT}}_\X$ acts
trivially. Let $D(-1) = D \otimes \Q(-1)$, a line in
$\Hom(\wedge^2V,\Q)$. Since ${\bigtilde{MT}}_\X$ acts trivially on $D$,
${\bigtilde{MT}}_\X$ acts on $D(-1)$ via $p_2$. Let
$$B = \{\alpha \in GL(V) : \alpha D(-1) \subseteq D(-1) \}$$
and let the character $\gamma : B \to \Aut(D(-1)) = {\mathbf G}_m$
be induced by the action of $GL(V)$ on $\Hom(\wedge^2V,\Q)$.
The action of $GL(V) \times {\mathbf G}_m$ on $\Hom(\wedge^2V,\Q)$ factors
through $GL(V)$. Therefore $MT_\X \subseteq B$, and we have a
commutative diagram
\downdiag{{\bigtilde{MT}}_\X}{MT_\X}{{\mathbf G}_m}{p_1}{p_2}{\gamma}
which gives the desired result.
\end{proof}

\begin{prop}
\label{rtilde}
If $\X$ is an abelian variety over a field $F$, $\ell$ is a prime number,
and $\ell \ne \fchar(F)$, then there exists a
(unique) character $\gamma_\ell : \r_\ell(F,\X) \to {\mathbf G}_m$ such that
\begin{enumerate}
\item[{(i)}]
${\tilde \r}_\ell(F,\X)$ is the graph of $\gamma_\ell$,
\item[{(ii)}]
the restriction of $\gamma_\ell$ to $G_{F,\X}$ is $\chi_\ell^{-1}$,
\item[{(iii)}]
if ${\bar F}$ is a subfield of $\C$, then $\gamma_\ell = \gamma$ on
$MT_{\X,\ell} \cap \r_\ell(F,\X)$.
\end{enumerate}
\end{prop}

\begin{proof}
Let $\pi_1$ and $\pi_2$ denote
the projection maps from $GL(V_\ell) \times {\mathbf G}_m$ onto $GL(V_\ell)$
and
${\mathbf G}_m$, respectively. By the definitions, $\r_\ell(F,\X)$ is the image
of ${\tilde \r}_\ell(F,\X)$ under $\pi_1$.
Fix a polarization on $\X$ defined over $F$. The
polarization generates a line $D_\ell$ in the $\Q_\ell$-vector space
$\Hom(\wedge^2V_\ell,\Q_\ell(1))$.
Let $D_\ell(-1) = D_\ell \otimes \Q_\ell(-1)$, a line in
$\Hom(\wedge^2V_\ell,\Q_\ell)$. Since the Weil pairing is
$\Gal(F^s/F)$-equivariant,
$\Gal(F^s/F)$ acts trivially on $D_\ell$. Therefore
${\tilde \r}_\ell(F,\X)$ acts trivially on $D_\ell$, and acts via $\pi_2$ on
$D_\ell(-1)$. Let
$$B_\ell =
\{\alpha \in GL(V_\ell) : \alpha D_\ell(-1) \subseteq D_\ell(-1) \}$$
and let the character
$\gamma_\ell : B_\ell \to \Aut(D_\ell(-1)) = {\mathbf G}_m$
be induced by the
action of $GL(V_\ell)$ on $\Hom(\wedge^2V_\ell,\Q_\ell)$.
The action of $GL(V_\ell) \times {\mathbf G}_m$ on
$\Hom(\wedge^2V_\ell,\Q_\ell)$ factors
through the action of $GL(V_\ell)$. Therefore
${\tilde \r}_\ell(F,\X) \subseteq B_\ell$, and we have a
commutative diagram
\downdiag{{\tilde \r}_\ell(F,\X)}{\r_\ell(F,\X)}{{\mathbf
G}_m}{\pi_1}{\pi_2}{\gamma_\ell}
which gives (i).
Since the restriction of $\pi_2$ to $G_{F,\X}$ is $\chi_\ell^{-1}$, we have
(ii). Now suppose ${\bar F}$ is a subfield of $\C$.
Using the fixed polarization, define $D$, $D(-1)$,
$B$, and $\gamma$ as in the proof of Theorem \ref{MTtilde}. Then
$B_\ell = B \times_\Q \Q_\ell$, and therefore
$MT_{\X,\ell} \subseteq B_\ell$. Since
$\gamma$ (respectively, $\gamma_\ell$) is induced by the action of
$GL(V)$ on $\Hom(\wedge^2V,\Q)$ (respectively, $GL(V_\ell)$ on
$\Hom(\wedge^2V_\ell,\Q_\ell)$), and $V_\ell = V \otimes_\Q \Q_\ell$,
we have (iii).
\end{proof}

Write ${\bigtilde{MT}}_{\X,\ell}$ for the $\Q_\ell$-algebraic subgroup
${\bigtilde{MT}}_\X \times_\Q \Q_\ell$ of $GL(V_\ell) \times {\mathbf G}_m$.
Then  ${\bigtilde{MT}}_{\X}(\Q_\ell) = {\bigtilde{MT}}_{\X,\ell}(\Q_\ell)$.
We state a reformulation of Theorem \ref{psdcor}, which we will use in
\S\ref{exclasses}.

\begin{thm}
\label{psdtilde}
If $\X$ is an abelian variety over a finitely generated extension $F$ of
$\Q$, then
${\tilde \r}_\ell(F,\X)^0 \subseteq {\bigtilde{MT}}_{\X,\ell}$.
\end{thm}

\begin{proof}
The result follows directly from Theorem \ref{psdcor} and Propositions
\ref{MTtilde} and \ref{rtilde}.
\end{proof}

\section{Connectedness extensions}
\label{exclasses}

\begin{thm}
\label{discond}
Suppose $\X$ is an abelian variety over a field $F$, $\ell$ is a prime
number not equal to $\fchar(F)$,
$$c : \Gal(F^s/F) \to \Aut_F(\X) \subseteq \Aut(V_\ell(\X))$$
is a homomorphism,
$\Y$ is the twist of $\X$ by the cocycle determined by $c$,
and $$F = F(\End(\X)) = F_{\Phi,\ell}(\X).$$
Then:
\begin{enumerate}
\item[{(i)}] $c$ induces an isomorphism
$$\Gal(F_{\Phi,\ell}(\Y)/F) \cong
\Im(c)/(\Im(c) \cap \r_\ell(F,\X)(\Q_\ell)),$$
\item[{(ii)}] $\r_\ell(F,\Y)$ is connected if and only
if $\Im(c) \subseteq \r_\ell(F,\X)(\Q_\ell)$,
\item[{(iii)}] if $M$ is the abelian extension of $F$ in
$F^s$ cut out by $c$, then $c$ induces an isomorphism
$$\Gal(M/F_{\Phi,\ell}(\Y)) \cong
\Im(c) \cap \r_\ell(F,\X)(\Q_\ell).$$
\end{enumerate}
\end{thm}

\begin{proof}
By Lemma \ref{conncomplem}ii, $F_{\Phi,\ell}(\Y) \subseteq M$.
The character $c$ induces isomorphisms
$$\Gal(M/F) \cong \Im(c)$$
and
$$\Gal(M/F_{\Phi,\ell}(\Y)) \cong
\Im(c) \cap \r_\ell(F,\Y)^0(\Q_\ell).$$
By Lemma \ref{conncomplem}i, we have
$\r_\ell(F,\Y)^0 \cong \r_\ell(F,\X)$, and the result follows.
\end{proof}

\begin{cor}
\label{discondcor}
Suppose $\X$ is an abelian variety over a field $F$, $\ell$ is a prime
number not equal to $\fchar(F)$, $\Y$ is the twist of $\X$ by a cocycle
$$c : \Gal(F^s/F) \to \Aut(\X) \subseteq \Aut(V_\ell(\X)),$$
and $$F = F(\End(\X)) = F_{\Phi,\ell}(\X) = F(\End(\Y)).$$
Then:
\begin{enumerate}
\item[{(i)}] $c$ is a character with values in $Z_\X^\times$
(where $Z_\X$ denotes the center of $\End(\X)$),
\item[{(ii)}] $c$ induces an isomorphism
$$\Gal(F_{\Phi,\ell}(\Y)/F) \cong
\Im(c)/(\Im(c) \cap \r_\ell(F,\X)(\Q_\ell)),$$
\item[{(iii)}] $\r_\ell(F,\Y)$ is connected if and only
if $\Im(c) \subseteq \r_\ell(F,\X)(\Q_\ell)$,
\item[{(iv)}] if $M$ is the abelian extension of $F$ in
$F^s$ cut out by $c$, then $c$ induces an isomorphism
$$\Gal(M/F_{\Phi,\ell}(\Y)) \cong
\Im(c) \cap \r_\ell(F,\X)(\Q_\ell).$$
\end{enumerate}
\end{cor}

\begin{proof}
By Lemma \ref{cocychar} and the assumption that
$F = F(\End(\Y))$, we have (i). The result now follows from
Theorem \ref{discond}.
\end{proof}

Suppose $A$ is an abelian variety defined over a field $F$ of characteristic
zero, $\k$ is a CM-field, $\iota : \k \hookrightarrow \End_F^0(A)$
is an embedding, and $C$ is an algebraically
closed field containing $F$. Let $\Lie(A)$ be the tangent space of $A$
at the origin,
an $F$-vector space. If $\sigma$ is an embedding of $\k$ into $C$, let
$$n_\sigma = \dim_C\{t \in \Lie(A)\otimes_F C :
\iota(\alpha)t = \sigma(\alpha)t {\text{ for all }} \alpha \in \k\}.$$
Write ${\bar \sigma}$ for the composition of $\sigma$ with the involution
complex conjugation of $\k$.

\begin{defn}
If $A$ is an abelian variety over an algebraically closed field $C$ of
characteristic zero, $\k$ is a CM-field, and
$\iota : \k \hookrightarrow \End^0(A)$ is an embedding,
we say $(A,\k,\iota)$ is {\em of Weil type}
if $n_\sigma = n_{\bar \sigma}$ for all embeddings $\sigma$ of $\k$ into
$C$.
\end{defn}

Although we do not use this fact, we remark that $(A,\k,\iota)$
is of Weil type if and only if $\iota$ makes $\Lie(A) \otimes_F C$ into a free
$\k \otimes_\Q C$-module (see p.~525 of \cite{Ribet} for the case where $\k$
is an imaginary quadratic field). Using the
semisimplicity of the $F$-algebra $\k \otimes_\Q F$ and the
$C$-algebra $\k \otimes_\Q C$, one may easily deduce that $\iota$ makes
$\Lie(A) \otimes_F C$ into a free $\k \otimes_\Q C$-module
if and only if $\iota$ makes $\Lie(A)$ into a free $\k \otimes_\Q F$-module.

Suppose $(A,\k,\iota)$ is of Weil type, and we have an element of
$$H^1(\Gal({\bar F}/F),\Aut(A))$$ which is represented by a cocycle $c$ with
values in the center of $\End^0(A)$. Let $B$ be the twist of $A$ by $c$,
and let $\varphi$ be the isomorphism from $\End_F(A)$ to $\End_F(B)$
obtained in Corollary \ref{cocy} and Proposition \ref{endaut}. Since
$(A,\iota)$ and
$(B,\varphi\circ \iota)$ are isomorphic over $C$, it follows that
$(B,\k,\varphi\circ \iota)$ is of Weil type.

Note that if $(A,\k,\iota)$ is of Weil type, then $\dim(A)$ is
divisible by $[\k:\Q]$.

\begin{thm}
\label{twistX}
Suppose $A$ is an abelian variety over a finitely generated
extension $F$ of $\Q$,
$\ell$ is a prime number, $\k$ is a CM-field,
$\iota$ is an embedding of $\k$ into the center of
$\End^0(A)$ such that $(A, \k, \iota)$ is of Weil type,
$c: \Gal({\bar F}/F) \to \k^\times$ is a character
of finite order $n$,
$r = 2\dim(A)/[\k:\Q] \in \Z$,
$M$ is the $\Z/n\Z$-extension of $F$ cut out by $c$, and
$B$ is the twist of $A$ by $c$.
Suppose
$F = F(\End(A))$, $\iota \circ c$ takes values in $\Aut(A)$,
$r$ is even, and $n$ does not divide $r$. Then
\begin{enumerate}
\item[{(i)}] $F = F(\End(B))$,
\item[{(ii)}]
either $F \ne F_{\Phi}(A)$ or $F \ne F_{\Phi}(B)$,
\item[{(iii)}] if $F_{\Phi}(A) = F$, then $F_{\Phi}(B) \subseteq M$
and $[M:F_{\Phi}(B)]$ divides $(n,2r)$,
\item[{(iv)}] if $F_{\Phi}(A) = F$ and $(n,2r) = 2$, then
$[M:F_{\Phi}(B)] = 2$.
\end{enumerate}
\end{thm}

\begin{proof}
The Galois module $V_{\ell}(B)$ is the twist of
$V_{\ell}(A)$ by $c$ (see Remark \ref{twistyrem}).
By applying Corollary \ref{cocy} to the cocycle induced
by $c$, we deduce (i) and we obtain an isomorphism $\varphi$
from $\End_F(A)$ onto $\End_F(B)$ such that $(B,\k,\varphi\circ\iota)$
is of Weil type.
Let $\k_\ell = \k \otimes \Q_\ell$. For $\CC = A$ or $B$, let
$$W_\CC =
\Hom_\Q(\wedge^{r}_{\k} H_1(\CC,\Q),\Q({\frac{r}{2}})), \quad
W_{\CC,\ell} =
\Hom_{\Q_\ell}(\wedge^{r}_{\k_\ell} V_\ell(\CC),\Q_\ell({\frac{r}{2}})),$$
where $\Hom_E$ means homomorphisms of $E$-vector spaces, if $E$ is a field.
Then $W_\CC$ is a one-dimensional $\k$-vector space and $W_{\CC,\ell}$ is a
free
rank-one
$\k_\ell$-module. The elements of $W_\CC$ are called Weil classes for $\CC$.
Since $V_{\ell}(\CC)=H_1(\CC,\Q)\otimes_{\Q}\Q_{\ell}$, we have
$W_{\CC,\ell}=W_\CC\otimes_{\Q}\Q_{\ell}.$
Consider the action of the Galois group $\Gal({\bar F}/F)$.
The Galois module
$W_{B,\ell}$ is the twist of the Galois module $W_{A,\ell}$
by the character $c^{-r}$.
Since $n$ does not divide $r$, this is a non-trivial twist,
so the Galois modules
$W_{B,\ell}$ and $W_{A,\ell}$ cannot be simultaneously trivial.

By pp.~52--54 of \cite{sln900} (see also Lemma 2.8 of \cite{moonzar}
and p.~423 of \cite{WeilHodge}), the elements of $W_\CC$ are Hodge classes
(since we are dealing with abelian varieties of Weil type).
Since ${\bigtilde MT}_\CC$ acts trivially on the Hodge classes,
${\bigtilde MT}_{\CC,\ell}(\Q_\ell)$ acts trivially on
$W_{\CC,\ell} = W_\CC \otimes_\Q \Q_\ell$.
Suppose now that $\r_\ell(F,A)$ and $\r_\ell(F,B)$ are both
connected. Then ${\tilde \r}_\ell(F,A)$ and ${\tilde \r}_\ell(F,B)$ are both
connected (by Proposition \ref{rtilde}i).
It follows from Theorem \ref{psdtilde} that
${\tilde \r}_\ell(F,\CC) \subseteq {\bigtilde MT}_{\CC,\ell}$.
Therefore, $W_{B,\ell}$ and $W_{A,\ell}$ are both trivial
as $\Gal({\bar F}/F)$-modules. This is a contradiction. We therefore have (ii).

Suppose that $F_{\Phi}(A) = F$. Then $\r_\ell(F,A)$ is connected,
so $\r_\ell(F,B)$ is disconnected. By Lemma \ref{conncomplem},
$\r_\ell(M,B)$ is connected. Therefore, $F_{\Phi}(B) \subseteq M$.
By Corollary \ref{discondcor}iv, $$\Gal(M/F_{\Phi}(B)) \cong
\Im(c) \cap \r_\ell(F,A)(\Q_\ell).$$
Let $\mu_s(k)$ denote the group of $s$-th roots of unity in $k^\times$.
We have $$\Im(c) = \mu_n(k) \cong \Z/n\Z.$$
Suppose $\alpha \in \Im(c) \cap \r_\ell(F,A)(\Q_\ell)$.
Then $\alpha^n = 1$. By Theorem \ref{psdcor} and the facts that
$\r_\ell(F,A) = \r_\ell(F,A)^0$ and
$\alpha \in \End^0(A)$, we have $\alpha \in MT_A(\Q)$.
Applying the character
$\gamma$ of Theorem \ref{MTtilde}, we have that $\gamma(\alpha)$
is an $n$-th root of unity in $\Q^\times$, and therefore
$\gamma(\alpha)$ is $1$ or $-1$.
By the definition of $W_{A,\ell}$, $\alpha$ acts on $W_{A,\ell}$
as multiplication by $\alpha^{-r}\gamma(\alpha)^{-r/2}$.
Since $\alpha \in MT_{A,\ell}(\Q_\ell)$,
$\alpha$ acts trivially on $W_{A,\ell}$. Therefore
$\alpha^{-r}\gamma(\alpha)^{-r/2} = 1$, so
$\alpha^{2r} = 1$.
Let $t= (n,2r)$.
Then
$$\Gal(M/F_{\Phi}(B)) \cong
\Im(c) \cap \r_\ell(F,A)(\Q_\ell) \subseteq
\mu_n(k) \cap \mu_{2r}(k) = \mu_t(k).$$
Therefore,
$[M:F_{\Phi}(B)]$ divides $t$.
Since $\r_\ell(F,\X)$ contains the homotheties
${\mathbf G}_m$ (see 2.3 of \cite{serrereps}), we have
$-1 \in \r_\ell(F,\X)(\Q_{\ell})$.
So $-1 \in \Im(c) \cap \r_\ell(F,\X)(\Q_{\ell})$
if and only if $-1 \in \Im(c)$, i.e., if and only if $n$ is even.
Thus if $t = 2$, then
$$\Gal(M/F_{\Phi}(B)) \cong \{\pm 1\}.$$
\end{proof}

\begin{thm}
\label{twistE}
Suppose $X$ and $Y$ are abelian varieties over a finitely generated
extension $F$ of $\Q$, $\ell$ is a prime number, $\Hom(X,Y) = 0$,
$F = F(\End(X)) =  F(\End(Y))$, $\k$ is a CM-field,
$[\k:\Q] = 2\dim(Y)$, and
$\dim(X) = t\dim(Y)$  for some odd positive integer $t$.
Suppose $\iota_X$
and $\iota_Y$ are embeddings of $\k$ into $\End^0(X)$ and $\End^0(Y)$,
respectively, and
$(X \times Y, \k, \iota_X \times \iota_Y)$
is of Weil type. Suppose $c$ is the non-trivial
character associated to a quadratic extension $M$ of $F$,
let $Y^c$ denote the twist of $Y$ by $c$, let $A = X \times Y$, and
let $B = X \times Y^c$. Then
\begin{enumerate}
\item[{(i)}] $F = F(\End(B))$,
\item[{(ii)}]
either $F(\End(A)) \ne F_{\Phi}(A)$ or $F(\End(B)) \ne F_{\Phi}(B)$,
\item[{(iii)}] if $F_{\Phi}(A) = F$, then $F_{\Phi}(B) = M$.
\end{enumerate}
\end{thm}

\begin{proof}
We have $F = F(\End(A))$.
Since $\Hom(X,Y) = 0$, we have $\End^0(A) = \End^0(X) \oplus \End^0(Y)$
and $\Aut(A) = \Aut(X) \times \Aut(Y)$. Consider the cocycle $c$ that
sends $\sigma \in \Gal({\bar F}/F)$ to $(1, c(\sigma)) \in
\Aut(X) \times \Aut(Y) = \Aut(A)$. All
the values of $c$ are of the form $(1, \pm 1)$, and therefore
belong to the center of
$\End^0(A)$. The abelian variety $B$ ( $ = X \times  Y^c$) is the twist
of $A$ ( $ = X \times Y$) by $c$. By Corollary \ref{cocy} we have
(i), and we obtain an isomorphism $\varphi$ from $\End_F(A)$ onto
$\End_F(B)$ such that $(B, \k, \varphi\circ(\iota_X \times \iota_Y))$ is of
Weil type.
Let $\k_\ell = \k \otimes \Q_\ell$.
We have
$$V_\ell(A) = V_\ell(X) \oplus V_\ell(Y), \quad
V_\ell(B) = V_\ell(X) \oplus V_\ell(Y^c).$$
Viewing the Tate modules as free $\k_\ell$-modules, we have
$$\wedge^{t+1}_{\k_{\ell}} V_\ell(A) = \wedge^{t}_{\k_\ell}V_\ell(X) \otimes
 _{\k_{\ell}} V_\ell(Y), \quad
\wedge^{t+1}_{\k_{\ell}} V_{\ell}(B)= \wedge^{t}_{\k_\ell}V_\ell(X)
\otimes_{\k_{\ell}} V_\ell(Y^c).$$  For
$\CC = A$ or $B$, let
$$W_{\CC,\ell} =
\Hom_{\Q_\ell}(\wedge^{t+1}_{\k_\ell} V_\ell(\CC),\Q_\ell({\frac{t+1}{2}})).$$
The Galois module $V_{\ell}(Y^{c})$ is the twist of the Galois module
$V_{\ell}(Y)$ by the character $c$ (see Remark \ref{twistyrem}),
and the Galois module
$W_{B,\ell}$ is the twist of the Galois module $W_{A,\ell}$ by
$c^{-1}=c$. Since $c$ is non-trivial, the Galois modules
$W_{B,\ell}$ and $W_{A,\ell}$ cannot be simultaneously trivial.
As in the proof of Theorem \ref{twistX}, it follows that
$\r_\ell(F,A)$ and $\r_\ell(F,B)$ cannot both be connected.
If $F_{\Phi}(A) = F$, then $\r_\ell(F,A)$ is connected, so
$\r_\ell(F,B)$ is disconnected. By Lemma \ref{conncomplem},
$\r_\ell(M,B)$ is connected, and so $F_{\Phi}(B)$ must be the
quadratic extension $M$ of $F$.
\end{proof}

\begin{rem}
Suppose $F$ is a subfield of $\C$, $Y$ is an elliptic curve over
$F$ with complex multiplication by an
imaginary quadratic field $K$, and $X$ is an absolutely simple
$3$-dimensional abelian
variety over $F$ with $K$ embedded in its endomorphism algebra.
Then we can always ensure (by taking complex
conjugates if necessary) that the two embeddings of
$K$ into $\C$ occur with the same multiplicity in the action of $K$ on the
tangent space of the $4$-dimensional abelian variety $A = X \times Y$.
Note that the hypotheses of Theorem \ref{twistE} (or of
Theorem \ref{twistX}) cannot be simultaneously
satisfied with $\dim(A) < 4$. In Example \ref{bigex} we exhibit
$4$-dimensional abelian
varieties satisfying the hypotheses of Theorem \ref{twistE}.
\end{rem}

\section{Examples}
\label{examples}

Using Theorems \ref{twistX} and \ref{twistE}, we can construct
examples of abelian varieties $B$ such that $\r_\ell(F,B)$
is disconnected, and compute the connectedness extensions.

\subsection{Example}
Let $\k = \Q(\sqrt{-3})$ and let $\kk$ be the CM-field which
is the compositum of
$\Q(\sqrt{-3})$ with the maximal totally real subfield $L$ of
$\Q(\zeta_{17})$. Then
$$\Gal({\kk}/\Q) \cong \Gal({\k}/\Q) \times \Gal({L}/\Q) \cong
\Z/2\Z \times \Z/8\Z.$$
Let $\Psi$ be the subset of $\Z/2\Z \times \Z/8\Z \cong
\Gal({\kk}/\Q)$ defined by
$$\Psi =
\{(0,0), (0,1), (0,4), (0,7), (1,2), (1,3), (1,5), (1,6)\}.$$
Let $\O_{\kk}$ denote the ring of integers of $\kk$. Let
$(A,\iota_{\kk})$ be an $8$-dimensional CM abelian variety of
CM-type $(\kk,\Psi)$
constructed from the lattice $\O_{\kk}$ as in Theorem 3
on p.~46 of \cite{st}, and defined over a number field $F$ (this can
be done by Proposition 26 on p.~109 of \cite{st}).
Then $A$ is absolutely simple, by the choice of $\Psi$ and
Proposition 26 on p.~69 of \cite{st}, and
$\End(A) = \O_{\kk}$ (see Proposition 6 on p.~42 of \cite{st}). Further,
the reflex field of $(\kk,\Psi)$ is $\kk$.
Take the number field $F$ to be
sufficiently large so that $F_\Phi(A) = F$. Let
$\iota$ be the restriction of $\iota_{\kk}$ to $\k$. By
the definition of $\Psi$, if
$\sigma \in \Gal(\k/\Q) = \Hom(\k,\C)$
then $n_\sigma = 4$. Therefore
$(A,\k,\iota)$ is of Weil type.
Let $c: \Gal({\bar F}/F) \to \k^\times$
be a non-trivial cubic character associated to a cubic extension
$M$ of $F$, and let $B$ denote
the twist of $A$ by $c$. Applying Theorem \ref{twistX}iii
with $n = 3$ and $r = 8$, then $F(\End(B)) = F$
and $F_\Phi(B) = M$.

\subsection{Example}
\label{bigex}
Let $J$ be the Jacobian of the genus $3$ curve $$y^7 = x(1-x),$$
and let $E$ be the elliptic curve $X_0(49)$. A model for $E$ is given
by the equation
$$y^2 + xy = x^3 - x^2 - 2x - 1.$$
Let $d$ be a non-zero square-free integer. If $d \ne 1$
let $E^{(d)}$ be the twist of $E$ by the non-trivial character of
$\Q(\sqrt{d})$, and if $d = 1$ let
$E^{(d)} = E$.
Let $$A = J \times E, \qquad A^{(d)} = J \times E^{(d)}.$$
The abelian varieties $A^{(d)}$ are defined over $\Q$. Let $\zeta_7$ be
a primitive seventh root of unity and let
$${\kk} = \Q(\zeta_7), \quad L_d = {\kk}(\sqrt{d}), \quad {\text{ and }} \quad
\k = \Q(\sqrt{-7}).$$
If $d = 1$ or $-7$ then ${\kk} = L_d$; otherwise, $[L_d:{\kk}] = 2$.
The abelian variety $J$ is a simple abelian variety with complex
multiplication by ${\kk}$, and the elliptic curves $E^{(d)}$ have complex
multiplication by the subfield $\k$ of ${\kk}$.
We have $\Gal({\kk}/\Q) = \{\sigma_1,...,\sigma_6\}$
where $\sigma_i(\zeta_7) = \zeta_7^i$.
The CM-type of $J$ is $({\kk}, \{\sigma_1, \sigma_2, \sigma_3\})$
(see p.~34 of \cite{Lang} or \S15.4.2 of \cite{st}), and the
reflex CM-type is $({\kk}, \{\sigma_4, \sigma_5, \sigma_6\})$
(see \S8.4.1 of \cite{st}).
We can identify $\End(A^{(d)})$ with the direct sum of
$\End(J)$ and $\End(E^{(d)})$. By Proposition 30 on p.~74 of \cite{st},
the smallest extension of $\Q$ over which all the elements of $\End(J)$ are
defined is the reflex CM-field of the CM-type of $J$, which is ${\kk}$.
Similarly, $ \k$ is
the smallest extension of $\Q$ over which all the elements of $\End(E^{(d)})$
are defined. We therefore have
$${\kk}  =  \Q(\End(A^{(d)})).$$

Next, we will prove that $L_d =  \Q_{\Phi}(A^{(d)})$.

Write $\O_\Omega$ for the ring of integers of a number field $\Omega$.
If $q$ is a prime number, let $\O_q = \O_\Omega \otimes \Z_q$.

\begin{claim}
\label{powerofp}
If ${\kk}'$ is a finite abelian extension of ${\kk}$ which is unramified away
from
the primes above $7$, then $[{\kk}':{\kk}]$ is a power of $7$.
\end{claim}

\begin{proof}
We have $-1 - \zeta_7 = (1 - \zeta_7^2)/(\zeta_7 - 1) \in \O_{\kk}^\times$.
Let $\P$ be the prime ideal of $\kk$ above $7$. The reduction map
$$\O_{\kk}^\times \to (\O_{\kk}/\P)^\times \cong (\Z/7\Z)^\times$$
is surjective, since $-1 - \zeta_7$ maps to $-2$, a generator of
$(\Z/7\Z)^\times$. Moreover, the class number of ${\kk}$ is one.
Therefore by class field theory, there is no
non-trivial abelian extension of ${\kk}$ of degree prime to $7$ and
unramified away from the primes above $\P$.
\end{proof}

\begin{claim}
\label{Kitself}
If $p$ is a prime and $p \equiv 3 \pmod{7}$,
then the only field ${\kk}'$ such that
\begin{enumerate}
\item[{(i)}] ${\kk}  \subseteq  {\kk}'  \subseteq  {\kk}(A_p)$,  and
\item[{(ii)}] ${\kk}'/{\kk}$ is unramified away from the primes above $7$,
\end{enumerate}
\noindent is ${\kk}$ itself.
\end{claim}

\begin{proof}
Since ${\kk}$ is a field of definition for the endomorphisms of the CM
abelian varieties $J$ and $E$, the extension ${\kk}(A_n)/{\kk}$ is abelian for
every integer $n$ (see Corollary 2 on p.~502 of \cite{serretate}).
Suppose $p$ and ${\kk}'$ satisfy the hypotheses of Lemma \ref{Kitself}.
Let $I_p \subseteq \Gal({\kk}(J_p)/{\kk})$ be the inertia subgroup at $p$.
We will first show
\begin{equation}
\label{claimb}
\#(I_p) = {\frac{p^6 - 1}{p^2 + p + 1}}.
\end{equation}
The image of $\O_p^\times$ in $\Gal({\kk}(J_p)/{\kk})$ under the Artin map of
class field theory is $I_p$, and we have natural homomorphisms
$$\Gal({\kk}(J_p)/{\kk}) \hookrightarrow
\Aut_{\O_{\kk}}(J_p) \cong
(\O_{\kk}/p\O_{\kk})^\times \cong \O_p^\times/(1 + p\O_p).$$
We therefore obtain maps
\begin{equation}
\label{maps}
\O_p^\times \to I_p \hookrightarrow \O_p^\times/(1 + p\O_p).
\end{equation}
Since the first map of (\ref{maps}) is surjective, the order of $I_p$ is
the order of the image of the composition.
Since $p \equiv 3 \pmod{7}$, we know that $p$ is inert in ${\kk}/\Q$, so
$(\O_{\kk}/p\O_{\kk})^\times$ ($\cong \O_p^\times/(1 + p\O_p)$) is a cyclic
group of
order $p^6 - 1$. Since the greatest common divisor of $p^6 - 1$ and
$p^3(p^2 + p + 1)$ is $p^2 + p + 1$, equation
(\ref{claimb}) will be proved when we show that the composition
of maps in (\ref{maps})
sends $u \in \O_p^\times$ to $u^{-p^3(p^2+p+1)} \pmod{1 + p\O_p}$.
We can view elements of $\Gal({\kk}/\Q)$ as automorphisms of $\O_p^\times$.
Proposition 7.40 on p.~211 of \cite{shimura} implies that the image of $u$
is of the form $\alpha(u)/\eta(u) \pmod{1 + p\O_p}$
where $\eta(u) = \sigma_4(u)\sigma_5(u)\sigma_6(u)$ and
$\alpha(u) \in {\kk}^\times$. Write ${\kk}_\A^\times$ for
the idele group of ${\kk}$, and for
each archimedean prime $\lambda$ of ${\kk}$, define
a Gr\"ossencharacter
$\psi_\lambda : {\kk}_\A^\times \to \C^\times$ by
$\psi_\lambda(x) =  (\alpha(x)/\eta(x))_\lambda$.
View $\O_p^\times$ as a subgroup of ${\kk}_\A^\times$.
Since $J$ has good reduction outside $7$, we have
$\psi_\lambda(\O_p^\times) = 1$,
by Theorem 7.42 of \cite{shimura}.
For $u \in \O_p^\times$, we have
$1 = \psi_\lambda(u) = \alpha(u)_\lambda = \alpha(u)$.
Therefore the image of $u$ in $\O_p^\times/(1 + p\O_p)$ is
$1/\eta(u) \pmod{1 + p\O_p}$. Since $p$ is inert in ${\kk}/\Q$, we have
$\Gal({\kk}/\Q) \cong \Gal((\O_{\kk}/p)/(\Z/p)) = D_p$, where $D_p$ is
the decomposition group at $p$. The latter group is a cyclic group of
order $6$ generated by the Frobenius element,
and we compute that
$$\sigma_4(u) \equiv u^{p^4}, \quad \sigma_5(u) \equiv u^{p^5},  \quad
{\text{and}}  \quad \sigma_6(u) \equiv u^{p^3} \quad \pmod{1 + p\O_p}$$
(since $p^4 \equiv 4 \pmod{7}$, $p^5 \equiv 5 \pmod{7}$,
and $p^3 \equiv 6 \pmod{7}$). Therefore
$$1/\eta(u) \equiv u^{-p^3(p^2+p+1)} \pmod{1 + p\O_p},$$
as desired.

We have
$$\Gal({\kk}(E_p)/{\kk}) \hookrightarrow
\Aut_{\O_\k}(E_p) \cong (\O_\k/p\O_\k)^{\times}.$$
The order of $(\O_\k/p\O_\k)^{\times}$ is $p^2 - 1$, which is not
divisible by $7$. Therefore $[{\kk}(A_p):{\kk}(J_p)]$ is not divisible by $7$.
By Lemma \ref{powerofp}, $[{\kk}':{\kk}]$ is a power of $7$. Therefore
${\kk}' \subseteq {\kk}(J_p)$. Since ${\kk}'/{\kk}$ is unramified
at $p$, we have
$I_p \subseteq \Gal({\kk}(J_p)/{\kk}')$.
Suppose ${\kk}' \neq {\kk}$. Then
$\#(I_p)$ divides $(p^6 - 1)/7$.  By (\ref{claimb}), $(p^6 - 1)/(p^2 + p + 1)$
divides $(p^6 - 1)/7$. Therefore $7$ divides $p^2 + p + 1$, which contradicts
the assumption that $p \equiv 3 \pmod{7}$. Therefore, ${\kk}' = {\kk}$.
\end{proof}

Suppose $p$ and $q$ are distinct odd primes, and $p \equiv 3 \pmod{7}$.
Let ${\kk}' = \kk(A_p) \cap \kk(A_q)$.
Since $A$ has good reduction outside $7$, the extension
${\kk}'/\kk$ is unramified away from the primes above $7$.
By Lemma \ref{Kitself}, we have ${\kk}' = \kk$.
As mentioned in the introduction, for
every integer $n \ge 3$ we have
$$\kk_{\Phi}(A) \subseteq \kk(A_n).$$
We therefore obtain
$$\kk = \kk_\Phi(A) = \Q_\Phi(A).$$
It follows from Theorem \ref{twistE} that
$$L_d =  \Q_{\Phi}(A^{(d)}).$$
Note that Shioda (see Theorem 4.4 of \cite{shioda}) proved the Hodge
Conjecture for $A$, and therefore
also for $A^{(d)}$. Thus, the Weil classes
on $A^{(d)}$ are algebraic. It follows easily that $L_d$ is the smallest
extension of $\Q$ over which all the algebraic cycle classes on all powers
of $A^{(d)}$ are defined.

\begin{remss}
If $\X$ is an abelian variety over a finitely generated extension $F$ of $\Q$,
and if the (as yet unproved) Tate Conjecture is true for all powers of $\X$
over $F_\Phi(\X)$,
then the field $F_\Phi(\X)$ is the smallest extension of $F$
over which all the algebraic cycle classes on all powers of $\X$ are
defined.
\end{remss}

\end{document}